\documentclass[12pt]{article}
\usepackage{e-jc}

\usepackage{amssymb}
\usepackage{amsmath, amsthm, enumitem}
\usepackage{graphicx}
\newtheorem{thm}{Theorem}

\newtheorem{claim}{Claim}
\newtheorem{prop}{Proposition}

\title{A Note on Total and Paired Domination of Cartesian Product Graphs}

\author{K. Choudhary\\
\small Department of Mathematics and Statistics\\[-0.8ex]
\small Indian Institute of Technology Kanpur\\[-0.8ex]
\small Kanpur, India\\
\small \texttt{keerti.india@gmail.com}\\
\and
S. Margulies\\
\small Department of Mathematics\\[-0.8ex]
\small Pennsylvania State University\\[-0.8ex]
\small State College, PA\\
\small \texttt{margulies@math.psu.edu}\\
\and
I. V. Hicks\\
\small Department of Computational and Applied Mathematics\\[-0.8ex]
\small Rice University\\[-0.8ex]
\small Houston, TX\\
\small \texttt{ivhicks@rice.edu}
}

\begin{document}

\maketitle

\begin{abstract}
A dominating set $D$ for a graph $G$ is a subset of $V(G)$ such that any
vertex not in $D$ has at least one neighbor in $D$.
The domination number $\gamma(G)$ is the size of a minimum dominating set in $G$. Vizing's conjecture from 1968 states that for the Cartesian product of graphs $G$ and $H$, $\gamma(G) \gamma(H) \leq \gamma(G \Box H)$, and Clark and Suen (2000) proved that $\gamma(G) \gamma(H) \leq 2\gamma(G \Box H)$. In this paper, we modify  the approach of Clark and Suen to prove a variety of similar bounds related to total and paired domination, and also extend these bounds to the $n$-Cartesian product of graphs $A^1$ through $A^n$.
\end{abstract}

%%%%%%%%%%%%%%%%%%%%%%%%%%%%%%%%%%%%%%%%%%%%%%%%%%%%%%%%%%%%%%%%%%%%%%%%
\section{Introduction}
We consider simple undirected graphs $G=(V,E)$ with vertex set $V$ and edge set $E$.
The open neighborhood of a vertex $v\in V(G)$ is denoted by $N_{G}(v)$,
and closed neighborhood by $N_{G}[v]$.
A dominating set $D$ of a graph G is a subset of $V(G)$ such that for all $v$,
$N_{G}[v] \cap D \neq \emptyset$.
A $\gamma$-set of G is a minimum dominating set for $G$, and its size is denoted
$\gamma(G)$.
A total dominating set $D$ of a graph $G$ is a subset of $V(G)$ such that for all $v$,
$N_{G}(v) \cap D \neq \emptyset$.
A $\gamma _{t}$-set of $G$ is a minimum total dominating set for $G$, and its size
is denoted $\gamma _{t}(G)$.
A paired dominating set $D$ for a graph $G$ is a dominating set such that the subgraph of $G$ induced by $D$ (denoted $G[D]$) has a perfect matching.
A $\gamma_{pr}$-set of $G$ is a minimum paired dominating set for $G$,
and its size is denoted $\gamma_{pr}(G)$.
In general, for a graph containing no isolated vertices,
$\gamma(G) \leq \gamma_{t}(G) \leq \gamma_{pr}(G)$. 
%For an easy example, consider %$P_5$, the path on five vertices, with vertices consecutively indexed from $1$ to %$5$. Then, a minimum dominating set consists of vertices $\{2,4\}$, a minimum 
%total dominating set consists of vertices $\{2,3,4\}$ and a minimum paired 
%dominating set consists of vertices $\{1,2,3,4\}$.

The Cartesian product graph, denoted $G \Box H$, is the graph with
vertex set $V(G)\times V(H)$, where vertices $gh$ and $g'h'$ are adjacent
whenever $g = g'$ and $(h,h')\in E(H)$, or $h = h'$ and $(g,g') \in E(G)$.
Just as the Cartesian product of graphs $G$ and $H$ is denoted $G \Box H$, the $n$-product of graphs $A^{1} , A^{2} ,\ldots, A^{n}$ is denoted as $A^{1} \Box A^{2} \Box \cdots \Box A^{n}$, and has vertex set
$V(A^{1})\times V(A^{2})\times \cdots \times V(A^{n})$, where vertices
$u^{1}\cdots u^{n}$ and $v^{1}\cdots v^{n}$ are adjacent if and only if
for some $i$, $(u^{i}, v^{i})\in E(A^{i})$, and $u^{j} = v^{j}$ for
all other indices $j \neq i$.

Vizing's conjecture from 1968 states that $\gamma(G) \gamma(H) \leq \gamma(G \Box H)$. For a thorough review of the activity on this famous open problem, see \cite{viz_survey_2009} and references therein. In 2000, Clark and Suen \cite{viz_clark_suen} proved that 
$\gamma(G) \gamma(H) \leq 2\gamma(G \Box H)$ by a sophisticated double-counting argument which involved projecting a $\gamma$-set of the product graph $G\Box H$ down onto the graph $H$. In this paper, we slightly modify the Clark and Suen double-counting approach and instead project subsets of $G \Box H$ down onto both graphs $G$ and $H$, which allow us to prove five theorems relating to total and paired domination. In this section, we state the results, and in Section \ref{sec_results}, we prove the results.

\begin{thm}~\label{thm1}
Given graphs $G$ and $H$ containing no isolated vertices,
\begin{align*}
max\big\{\gamma(G) \gamma _{t}(H), \gamma _{t}(G) \gamma(H)\big\} \leq 2\gamma(G \Box H)~.
\end{align*}
\end{thm}

In 2008, Ho \cite{ho_total_dom} proved an inequality for total
domination analogous to the Clark and Suen inequality for domination. In particular, Ho proved $\gamma_{t}(G) \gamma _{t}(H) \leq 2\gamma_{t}(G \Box H)$. We provide a slightly different proof of Ho's inequality, and then extend the result to the $n$-product case.

\begin{thm}[Ho~\cite{ho_total_dom}]~\label{thm2}
Given graphs $G$ and $H$ containing no isolated vertices,
\begin{align*}
\gamma_{t}(G) \gamma _{t}(H) \leq 2\gamma_{t}(G \Box H)~.
\end{align*}
\end{thm}

\begin{thm}~\label{thm3}
Given graphs $A^{1} , A^{2} ,\ldots, A^{n}$ containing no isolated vertices,
\begin{align*}
\displaystyle \prod\limits _{i=1}^{n} \gamma_{t}(A^{i}) \leq n\gamma_{t}(A^{1} \Box A^{2} \Box \cdots \Box A^{n})~.
\end{align*}
\end{thm}
In 2010, Hou and Jiang \cite{hou_jiang_pr_dom} proved that $\gamma_{pr}(G) \gamma _{pr}(H) \leq 7\gamma_{pr}(G \Box H)$, for graphs $G$ and $H$ containing no isolated vertices. We provide an improvement to this result, and extend the result to the $n$-product graph.
\begin{thm}~\label{thm4}
Given graphs $G$ and $H$ containing no isolated vertices,
\begin{align*}
\gamma_{pr}(G) \gamma _{pr}(H) \leq 6 \gamma_{pr }(G \Box H)~.
\end{align*}
\end{thm}

\begin{thm}~\label{thm_pr_n}
Given graphs $A^{1},\ldots, A^{n}$ containing no isolated vertices,
\begin{align*}
\displaystyle\prod\limits_{i=1}^{n} \gamma _{pr}(A_i) \leq 2^{n-1}(2n - 1)\gamma _{pr}(A_1 \Box \cdots \Box A_n)~.
\end{align*}
\end{thm}

%%%%%%%%%%%%%%%%%%%%%%%%%%%%%%%%%%%%%%%%%%%%%%%%%%%%%%%%%%%%%%%%%%%%%%%%
\section{Main Results} \label{sec_results}
We begin by introducing some notation which will be utilized throughout the proofs in this section. Given $S \subseteq V(G\Box H)$, the projection of $S$ onto graphs $G$ and $H$ is defined as
\begin{align*}
\varPhi_{G}(S) &= \{g \in V(G) ~\rvert~\exists ~h \in V(H) \text{ with } gh \in S\}~,\\
\varPhi_{H}(S) &= \{h \in V(H) ~\rvert~\exists ~g \in V(G) \text{ with } gh \in S\}~.
\end{align*}
In the case of the $n$-product graph $A^1\Box \cdots \Box A^n$, we project a set of vertices in $V(A^1\Box \cdots \Box A^n)$ down to a particular graph $A_i$. Therefore, given $S \subseteq V(A^1\Box \cdots \Box A^n)$, we define
\begin{align*}
\varPhi_{A^i}(S) &= \{a \in V(A^i) ~\rvert~\exists ~u^1\cdots u^n \in S \text{ with } a = u^i\}~.
\end{align*}
For $gh \in V(G\Box H)$, the $G$-neighborhood and $H$-neighborhood of $gh$ are defined as follows:
\begin{align*}
N_{\underline{\mathbf{G}} \Box H}(gh) &= \{g'h \in V(G\Box H) ~\rvert~ g' \in N_G(g)\}~,\\ 
N_{G \Box \underline{\mathbf{H}}}(gh) &= \{gh' \in V(G\Box H) ~\rvert~ h' \in N_H(h)\}~. 
\end{align*}
Thus, $N_{\underline{\mathbf{G}} \Box H}(gh)$ and $N_{G \Box \underline{\mathbf{H}}}(gh)$ are both subsets of $V(G \Box H)$. Additionally, $E(G \Box H)$ can be partitioned into two sets, \textbf{G}-edges and \textbf{H}-edges, where
\begin{align*}
\text{\textbf{G}-edges} &= \{(gh,g'h) \in E(G\Box H) ~\rvert~ h \in V(H) \text{ and } ( g,g')\in E(G)\}~, \\
\text{\textbf{H}-edges} &= \{(gh,gh')\in E(G\Box H) ~\rvert~ g \in V(G) \text{ and } (h,h')\in E(H)\}~.
\end{align*}
In the case of the $n$-product graph $A^1\Box \cdots \Box A^n$, we identify the $i$-neighborhood of a particular vertex, and partition the set of edges $E(A^1\Box \cdots \Box A^n)$ into $n$ sets. Thus, we define $E_i$ to be
\begin{align*}
E_i &= \Big\{\big(u^{1}\cdots u^{n}, v^{1}\cdots v^{n}\big) ~\lvert~ (u^{i}, v^{i})\in E(A^{i}) \text{, and $u_{j} = v_{j}$, for all other indices $j\neq i$}\Big\}~,
\end{align*}
and for a vertex $u \in V(A^1\Box \cdots \Box A^n)$, we define
\begin{align*}
N_{\Box A^i}(u) &= \Big\{ v\in  V(A^1\Box \cdots \Box A^n) ~\lvert~\text{$v$ and $u$ are connected by $E_{i}$-edge} \Big\}~.
\end{align*}
Finally, we need two elementary propositions about matrices that will be utilized throughout the proofs. 
\begin{prop}\label{fact1}
Let $M$ be a binary matrix. Then either
%\begin{enumerate}[label=(\roman*)]
\begin{enumerate}[label=(\alph*)]
\item each column contains a 1, or 
\item each row contains a 0~.
\end{enumerate}
\end{prop}
Prop. 1 refers only to $d_1 \times d_2$ binary matrices. Prop. 2 is a generalization of Prop. 1 for $d_1 \times d_2 \times \cdots \times d_n$ $n$-ary matrices. 
\begin{prop}\label{fact2}
Let $M$ be a $d_1 \times d_2 \times \cdots \times d_n$, $n$-ary matrix ($n$-ary in this case signifies that $M$ contains entries only in the range $\{1,\ldots, n\}$). Then there exists a $j \in \{1,\ldots, n\}$ (not necessarily unique),
such that each of the $d_1 \times \cdots \times d_{j-1} \times 1 \times d_{j+1} \times \cdots \times d_n$ submatrices of $M$ contains an entry with value $j$. Such a matrix $M$ is called a $j$-matrix.
\end{prop}

Note that, given any $d_1 \times d_2 \times \cdots \times d_n$ matrix, there are $d_j$ submatrices of the form $d_1 \times \cdots \times d_{j-1} \times 1 \times d_{j+1} \times \cdots \times d_n$. We will denote such a submatrix as $M[:,i_j,:]$ with $1 \leq i_j \leq d_j$.

\begin{proof}
Let $M$ be a $d_1 \times d_2 \times \cdots \times d_n$ $n$-ary matrix which is not a $j$-matrix for $1 \leq j \leq n-1 $. We will show that $M$ is an $n$-matrix.

%Consider $j = 1$. Since $M$ is not a 1-matrix, each of the $d_2 \times d_3 \times \cdots \times d_n$ submatrices does \emph{not} contain a 1. Note that there are $d_1$ such submatrices. Next, consider $j = 2$. Since $M$ is also not a 2-matrix, each of the $d_1 \times d_3 \times \cdots \times d_n$ submatrices does \emph{not} contain a 2. Note that there are $d_2$ such submatrices. We continue this pattern for $1 \leq j \leq n - 1$. Since $M$ is \emph{not} a $j$-matrix for $1 \leq j \leq n-1$, $M$ does \emph{not} contain any elements in the set $\{1,\ldots, n-1\}$. But $M$ is an $n$-ary matrix. Thus, the $d_1 \times d_2 \times \cdots \times d_n$ matrix consists entirely of entries with value $n$. Thus, $M$ is an $n$-matrix.
%
%this is not true - Thus, the $d_1 \times d_2 \times \cdots \times d_n$ matrix consists entirely of entries with value $n$. 
%
Consider $j = 1$. Since $M$ is not a 1-matrix,
there exists at least one $1 \times d_2 \times d_3 \times \cdots \times d_n$
submatrix that does \emph{not} contain a 1. Without loss of generality, let
 $M[i_1,:]$ with $1 \leq i_1 \leq d_1$ be such a matrix.
Next, consider $j = 2$. Since $M$ is also not a 2-matrix, let 
$M[:,i_2,:]$ with $1 \leq i_2 \leq d_2$ be a $d_1 \times 1 \times d_3 \times \cdots \times d_n$ submatrix that does \emph{not} contain a 2.
Therefore, $M[i_1,i_2,:]$
is a $1 \times 1 \times d_3 \times \cdots \times d_n$
submatrix that contains neither a 1 nor a 2.
We continue this pattern for $1 \leq j \leq n - 1$.
Since $M$ is \emph{not} a $j$-matrix for $1 \leq j \leq n-1$, let
$M[i_1,\ldots,i_{n-1},:]$ be the $1 \times \cdots 1 \times d_n$ submatrix containing no elements in the set $\{1,\cdots, n-1\}$. 
Therefore, for all $1 \leq x \leq d_n$, $M[i_1, \ldots, i_{n-1}, x] = n$, and \emph{all} of the $d_1 \times \cdots \times d_{n-1} \times 1$ submatrices of $M$ contains an entry with value $n$. Thus, $M$ is an $n$-matrix.
\end{proof}

Now, we present the proofs of Theorems 1 through 5.

%%%%%%%%%%%%%%%%%%%%%%%%%%%%%%%%%%%%%%%%%%%%%%%%%%%%%%%%%%%%%%%%%%%%%%%%
\subsection{Proof of Theorem~\ref{thm1}}
\begin{proof}
Let ${\{u_{1} ,\ldots, u _{\gamma_{t}(G)}\} }$ be a $\gamma _{t}$-set of $G$.
%partitioning is slightly different for G and H.
Partition $V(G)$ into sets $D _{1} ,\ldots, D _{\gamma_{t}(G)}$,
such that $D _{i} \subseteq N_{G}(u _{i})$. Let ${\{\overline{u} _{1} ,\ldots, \overline{u} _{\gamma(H)}\} }$ be a $\gamma$-set of $H$.
Partition $V(H)$ into sets $\overline{D} _{1} ,\ldots, \overline{D} _{\gamma(H)}$,
such that $\overline{u} _{j} \in \overline{D} _{j}$ and $\overline{D} _{j} \subseteq N_{H}[\overline{u} _{j}]$. We note that ${\{D _{1} ,\ldots, D _{\gamma _{t}(G)}\} } \times {\{\overline{D} _{1} ,\ldots, \overline{D} _{\gamma(H)}\} }$ is a partition of $V(G\Box H)$. Let $D$ be a $\gamma$-set of $G \Box H$. Then, for each $gh \notin D$, either $N_{\underline{\mathbf{G}} \Box H}(gh) \cap D$
or $N_{G \Box \underline{\mathbf{H}}}(gh) \cap D$ is non-empty. Based on this observation, we define the binary $|V(G)| \times |V(H)|$ matrix $F$ such that:
\[
 F(g,h) =
 \begin{cases}
 1 & \text{if } gh\in D\text{ or } N_{G \Box \underline{\mathbf{H}}}(gh) \cap D \neq \emptyset~, \\
 0 & \text{otherwise}~.
 \end{cases}
\]
Since $F$ is a $|V(G)| \times |V(H)|$ matrix, each of the $D_i \times \overline{D}_j$ subsets of $V(G \Box H)$ determines a submatrix of $F$.

For $i=1,\ldots,\gamma_{t}(G)$, let $Z _{i} = D \cap (D _{i} \times V(H))$, and let
\begin{align*}
S _{i} = \big\{\overline{D} _{x} ~\rvert~ &\text{the submatrix of $F$ determined by $D _{i}\times \overline{D} _{x}$ satisfies Prop.~\ref{fact1}a},\\
&\text{with $x \in \{1,\ldots, \gamma(H)\}$} \big\}~.
\end{align*}

For $j=1,\ldots,\gamma(H)$, let $\overline{Z} _{j} = D \cap (V(G) \times \overline{D} _{j})$, and let
\begin{align*}
\overline{S} _{j} = \big\{D _{x} ~\rvert~ &\text{the submatrix of $F$ determined by $D _{x}\times \overline{D} _{j}$ satisfies Prop.~\ref{fact1}b},\\
&\text{with $x \in \{1,\ldots, \gamma_t(G)\}$}\big\}~.
\end{align*} 

Let $d_H = \sum_{i=1}^{\gamma _{t}(G)} |S _{i} |$, and
$d_G = \sum_{j=1}^{\gamma(H)}  |\overline{S} _{j} |$. Since the partition of $V(G \Box H)$ composed of elements $D _{i}\times \overline{D} _{j}$ contains $\gamma _{t}(G) \gamma(H)$ components, and since every $D _{i}\times \overline{D} _{j}$ submatrix of $F$ satisfies either conditions (a) or (b) of Prop.~\ref{fact1} (possibly both), $\gamma _{t}(G) \gamma(H) \leq d_H + d_G $. We will now prove two subclaims which will allow us to bound the size of our various sets.

\begin{claim}\label{clm_dj}
If the submatrix of $F$ determined by $D _{i} \times \overline{D} _{j}$ satisfies Prop.~\ref{fact1}a, then
$\overline{D} _{j}$ is dominated by $\varPhi_{H}(Z _{i})$.
\end{claim}

\begin{proof}
Let $h \in \overline{D} _{j}$. We must show that either $h \in \varPhi_{H}(Z _{i})$, or $h$ is adjacent to a vertex $h'$ in $\varPhi_{H}(Z _{i})$. If $(D _{i}\times \{h\})\cap D \neq \emptyset$, there exists a $g \in D_i$ such that $gh \in D$. Thus, $h \in \varPhi_{H}(Z _{i})$. 

If $(D _{i}\times \{h\})\cap D = \emptyset$, then recall that the submatrix of $F$ determined by $D_i \times \overline{D} _{j}$ satisfies Prop.~\ref{fact1}a. Therefore, there is a 1 in every column of the submatrix. This implies there exists a $g \in D_i$ such that $F(g,h) = 1$. Since $gh \notin D$, there exists an $h' \in V(H)$ such that $gh' \in N_{G \Box \underline{\mathbf{H}}}(gh) \cap D$. Therefore, $(gh', gh)$ is an \textbf{H}-edge, implying $(h,h') \in E(H)$ and $h$ is adjacent to $h'$. Therefore, $\overline{D} _{j}$ is dominated by $\varPhi_{H}(Z _{i})$.
\end{proof}

\begin{claim}\label{clm_di_out}
If the submatrix of $F$ determined by $D _{i} \times \overline{D} _{j}$ satisfies Prop.~\ref{fact1}b, then
$D _{i}$ is dominated by $\varPhi_{G}(\overline{Z} _{j})$. Additionally, $\forall g \in D_i \cap \varPhi_{G}(\overline{Z}_{j})$, there exists a vertex $g' \in \varPhi_{G}(\overline{Z}_{j})$ such that $(g,g') \in E(G)$.
\end{claim}

We note that this claim does not imply that $\varPhi_{G}(\overline{Z} _{j})$ is a total dominating set, but the claim is a slightly stronger condition on domination. When applying this condition, we will say that the set $D_i$ is \emph{non-self dominated} by $\varPhi_{G}(\overline{Z}_{j})$.

\begin{proof} The argument for proving that $\varPhi_{G}(\overline{Z} _{j})$ dominates $D_i$ is almost identical to the proof of Claim \ref{clm_dj}. The only difference is that the $D _{i} \times \overline{D} _{j}$ submatrix of $F$ satisfies Prop.~\ref{fact1}b. Thus, every row contains a 0. But since every vertex in $V(G \Box H)$ is dominated by $D$, this implies that every vertex $g \in D_i$ is dominated by some other (not itself) vertex $g' \in \varPhi_{G}(\overline{Z} _{j})$. Thus, $D_i$ is dominated by $\varPhi_{G}(\overline{Z} _{j})$, with the slightly stronger condition that \emph{every} vertex in $D_i$ \big(even those vertices in $D_i \cap \varPhi_{G}(\overline{Z}_{j})$\big) is adjacent to \emph{another} vertex in $\varPhi_{G}(\overline{Z}_{j})$.
\end{proof}

\begin{claim}\label{claim_zi_si}
For $i=1,\ldots,\gamma _{t}(G)$, $|S_i| \leq | Z _{i}|$.
Similarly, for $j=1,\ldots,\gamma(H)$,
$|\overline{S} _{j}| \leq  |\overline{Z} _{j}|$.
\end{claim}
\begin{proof}
Let $S _{i} = \{\overline{D} _{j _{1}}, \overline{D} _{j _{2}}, \ldots, \overline{D} _{j _{k}}\}$,
and let $A = \varPhi_{H}(Z _{i})$. Note that $|A| \leq |Z_i|$.
By Claim \ref{clm_dj}, $A$ dominates $\cup _{x=1}^{k}\overline{D} _{j _{x}}$. Therefore,
$A \cup \big\{ \overline{u} _{j}~\rvert~ j\notin \{j _{1}, j _{2}, \ldots, j _{k}\}\big\}$
is a dominating set of $H$, and, since the sets $A$ and
$\big\{ \overline{u} _{j} ~\rvert~ j\notin \{j _{1}, j _{2}, \ldots, j _{k}\} \big\}$ are disjoint, then
\begin{align*}
\big\lvert A \cup \big\{ \overline{u} _{j} ~\rvert~ j\notin \{j _{1}, j _{2}, 
\ldots, j _{k}\}\big\} \big\lvert = \lvert A \lvert + ( \gamma(H) - k ) \geq \gamma(H)~.
\end{align*}
Hence, $ k = |S_i| \leq |A| \leq |Z _{i}|$.

For the proof of second part, let
$\overline{S} _{j}$ = $\{D _{i _{1}}, D _{i_{2}}, \ldots, D _{i _{k}}\}$,
and let $A$ = $\varPhi_{G}(\overline{Z} _{j})$.
Again, note that $ |A| \leq |\overline{Z} _{j}|$.
Then by Claim~\ref{clm_di_out}, $A$ dominates $\cup _{x=1}^{k} D _{i_{x}}$,
with the stronger condition that $\forall g \in D_{i_x} \cap A$,
there exists a vertex $g' \in A$ such that $(g,g') \in E(G)$.
Now we consider $A \cap \big\{ u _{i}~\rvert~i\notin \{i _{1}, i _{2}, \ldots, i _{k}\}\big\}$. If this intersection is non-empty, let $A \cap \big\{ u _{i}~\rvert~i\notin \{i _{1}, i _{2}, \ldots, i _{k}\}\big\} = \{u_{i_{k+1}},\ldots, u_{i_{l}}\}$. Then, $A$ dominates $\cup _{x=1}^{l} D _{i_{x}}$ with the same stronger condition. Moreover, the sets $A$ and $\big\{ u _{i}~\rvert~i\notin \{i _{1}, i _{2}, \ldots, i _{k}, \ldots, i _{l}\}\big\}$ are disjoint.

We claim that
$A \cup \big\{u _{i}~\rvert~i\notin \{i _{1}, \ldots, i _{l}\}\big\}$
is a total dominating set of $G$.
To see this, consider any vertex $g \in V(G)$.
If $g \in D_x$ with $x \in \{i _{1}, i _{2}, \ldots, i_{k}\}$,
then by the stronger condition on domination associated with Claim \ref{clm_di_out},
$g$ is adjacent to another vertex in $A$.
If $g \in D_x$ with $x \notin \{i _{1}, \ldots, i _{k}\}$,
then $u_x \in \big\{u _{i}~\rvert~i\notin \{i _{1}, \ldots, i_{k}\}\big\}$,
and $g$ is adjacent to $u_x$, since $u_x$ dominates $D_x$. We note that $u_x$ is either in $A$ (if $k+1\leq x \leq l$) or in $\big\{u _{i}~\rvert~i\notin \{i _{1}, \ldots, i_{l}\}\big\}$. In either case, $A \cup \big\{u _{i}~\rvert~i\notin \{i _{1}, \ldots, i_{l}\}\big\}$
is a total dominating set of $G$, and
\begin{align*}
\big \lvert A\cup \big\{ u _{i}~\rvert~i\notin \{i _{1}, i_{2}, \ldots, i _{l}\}\big\} \big\lvert = \lvert A \lvert + ( \gamma _{t}(G) - l ) \geq\gamma _{t}(G)~.
\end{align*}
Hence, as before, $ k = |\overline{S} _{j}| \leq l \leq |A| \leq | \overline{Z} _{j}|$~.
\end{proof}

To conclude the proof, we observe that
\begin{align*} 
d_H &= \displaystyle\sum\limits_{i=1}^{\gamma _{t}(G)} | S _{i} | \leq \displaystyle\sum\limits_{i=1}^{\gamma _{t}(G)} | Z _{i}| \leq | D|~,\\\
d_G &= \displaystyle\sum\limits_{j=1}^{\gamma(H)}  | \overline{S} _{j}| \leq \displaystyle\sum\limits_{j=1}^{\gamma(H)}  | \overline{Z} _{j} | \leq | D|~.
\end{align*}
Hence,
$\gamma _{t}(G) \gamma(H) \leq d_H
+ d_G \leq 2|D| \leq 2\gamma(G \Box H)$. Moreover, we can similarly prove that $\gamma(G) \gamma_{t}(H) \leq 2\gamma(G \Box H)$.
Therefore, $\max\{\gamma(G) \gamma _{t}(H), \gamma _{t}(G) \gamma(H)\} \leq 2\gamma(G \Box H)$.
\end{proof}

%%%%%%%%%%%%%%%%%%%%%%%%%%%%%%%%%%%%%%%%%%%%%%%%%%%%%%%%%%%%%%%%%%%%%%%%
\subsection{Proof of Theorem~\ref{thm2}}
\begin{proof}
Let ${\{u _{1} ,\ldots, u _{\gamma_{t}(G)}\} }$ be a $\gamma _{t}$-set of $G$.
Partition $V(G)$ into sets $D _{1} ,\ldots, D _{\gamma_{t}(G)}$,
such that if $u \in D _{i}$ then $u\in N_{G}(u _{i})$
for all $i = 1, \ldots, \gamma_{t}(G)$. Similarly,
let ${\{\overline{u} _{1} ,\ldots, \overline{u} _{\gamma_{t}(H)}\} }$ be a
$\gamma_{t}$-set of $H$ and $\overline{D} _{1} ,\ldots, \overline{D} _{\gamma_{t}(H)}$
be the corresponding partitions. Then,
${\{D _{1} ,\ldots, D _{\gamma _{t}(G)}\} } \times {\{\overline{D} _{1} ,\ldots, \overline{D} _{\gamma_{t}(H)}\} }$
forms a partition of $V(G\Box H)$.

Let $D$ be a $\gamma_{t}$-set of $G \Box H$. Then, for each $gh \in V(G \Box H)$, either the set $N_{\underline{\mathbf{G}} \Box H}(gh) \cap D$
or the set $N_{G \Box \underline{\mathbf{H}}}(gh) \cap D$ is non-empty. Based on this observation, we define the binary $|V(G)| \times |V(H)|$ matrix $F$:
\[
 F(g,h) =
 \begin{cases}
 1 & \text{if } N_{G \Box \underline{\mathbf{H}}}(gh) \cap D \neq \emptyset~, \\
 0 & \text{otherwise}~.
 \end{cases}
\]
For $i=1,\ldots,\gamma_{t}(G)$, let $Z _{i} = D \cap (D _{i} \times V(H))$, and let
\begin{align*}
S _{i} = \big\{\overline{D} _{x} ~\rvert~ &\text{the submatrix of $F$ determined by $D _{i}\times \overline{D} _{x}$ satisfies Prop.~\ref{fact1}a},\\
&\text{with $x \in \{1,\ldots, \gamma_t(H)\}$} \big\}~.
\end{align*}
For $j=1,\ldots,\gamma_{t}(H)$, let $\overline{Z} _{j} = D \cap (V(G) \times \overline{D} _{j})$, and let
\begin{align*}
\overline{S} _{j} = \big\{D _{x} ~\rvert~ &\text{the submatrix of $F$ determined by $D _{x}\times \overline{D} _{j}$ satisfies Prop.~\ref{fact1}b},\\
&\text{with $x \in \{1,\ldots, \gamma_t(G)\}$}\big\}~.
\end{align*}

Let $d_H = \sum_{i=1}^{\gamma _{t}(G)} |S _{i} |$, and
$d_G = \sum_{j=1}^{\gamma_t(H)}  |\overline{S} _{j} |$. Since the partition of $V(G \Box H)$ composed of elements $D _{i}\times \overline{D} _{j}$ contains $\gamma _{t}(G) \gamma_t(H)$ components, and since every submatrix of $F$ determined by $D _{i}\times \overline{D} _{j}$ satisfies either Prop.~\ref{fact1}a or \ref{fact1}b (or possibly both), then $\gamma _{t}(G) \gamma_t(H) \leq d_H + d_G $. 

Furthermore, by similar arguments given in the proof of Theorem~\ref{thm1} (specifically, Claims \ref{clm_dj} and \ref{clm_di_out}), we can conclude, as before, that 
for $i=1,\ldots,\gamma _{t}(G)$, $| S _{i}| \leq |Z _{i}|$ and, for $j=1,\ldots,\gamma_{t}(H)$,
$|\overline{S} _{j}| \leq |\overline{Z} _{j}|$. Finally,

\begin{align*}
d_H &= \sum_{i=1}^{\gamma _{t}(G)} |S _{i}| \leq   \sum_{i=1}^{\gamma _{t}(G)} | Z _{i}| = |D| = \gamma_t(G\Box H)~,\\
d_G &= \sum_{j=1}^{\gamma _{t}(H)} |\overline{S} _{j}| \leq   \sum_{j=1}^{\gamma _{t}(H)} | \overline{Z} _{j}| = |D| = \gamma_t(G\Box H)~.
\end{align*}
Summing these two equations, we see $d_H + d_G \leq 2 \gamma_t(G\Box H)$, which implies $\gamma _{t}(G) \gamma_{t}(H) \leq 2 \gamma_t(G\Box H)$~.
\end{proof}

%%%%%%%%%%%%%%%%%%%%%%%%%%%%%%%%%%%%%%%%%%%%%%%%%%%%%%%%%%%%%%%%%%%%%%%%
\subsection{Proof of Theorem~\ref{thm3}}
\begin{proof}
%Referring to the definition of $E_i$ given by Eq. \ref{eq_ei}, the set $\{E_1, %E_2,\ldots, E_n \}$ forms a partition of $E(A^{1} \Box \cdots \Box A^{n})$.
For $i=1,\ldots,n$, let ${\{u^{i} _{1} ,..., u^{i} _{\gamma_{t}(A^{i})}\} }$ be a
$\gamma _{t}$-set of $A^{i}$, and $D^{i} _{1} ,\ldots, D^{i} _{\gamma_{t}(A^{i})}$
be the corresponding partitions (as defined in the proof of Theorem \ref{thm2}).

Let $Q = \{D^{1} _{1} ,\ldots, D^{1} _{\gamma _{t}(A^{1})}\} \times \cdots \times \{D^{n} _{1} ,\ldots, D^{n} _{\gamma _{t}(A^{n})}\}$.
Then $Q$ forms a partition of $V(A^{1} \Box \cdots \Box A^{n})$ with 
$|Q| = \displaystyle\prod\limits_{i=1}^{n} \gamma_{t}(A^{i})$. 

Let $D$ be a $\gamma_{t}$-set of $A^{1} \Box \cdots \Box A^{n}$. Then, for each $u\in V(A^{1} \Box \cdots \Box A^{n})$, there exists an $i$ such that $N _{\Box A^i}(u) \cap D$ is non-empty. Based on this observation (as in the 2-dimensional case), we define an $n$-ary $|V(A^1)| \times \cdots \times |V(A^n)|$ matrix $F$ such that:
\begin{align*}
F(u_1,\dots,u_n) &= \min \{i~\lvert~N _{\Box A^i}(u_1\cdots u_n) \cap D \neq \emptyset \}~.
\end{align*}
For $j=1,\ldots,n$, let $d_{j} \subseteq Q$ be the set of the elements in $Q$ which are $j$-matrices. By Prop.~\ref{fact2}, each element of $Q$ belongs to at least one $d_{j}$-set. Then, $\displaystyle\prod\limits_{i=1}^{n} \gamma_{t}(A^{i}) \leq \displaystyle\sum\limits_{j=1}^{n} | d_{j}|$.

\begin{claim}\label{claim3}
For $j=1,\ldots,n$, $ |d_{j}| \leq |D|$.
\end{claim}
\begin{proof}
We prove here that $ |d_{n}| \leq |D|$, but a similar proof can be performed for any other $j$. Similar to $Q$, let $B = \{D^{1} _{1}, \ldots, D^{1} _{\gamma _{t}(A^{1})}\} \times \cdots \times \{D^{n-1} _{1} ,\ldots, D^{n-1} _{\gamma _{t}(A^{n-1})}\}$.
For convenience, we denote $B$ as $\{B_{1}, \ldots, B_{|B|}\}$, where $\rvert B\rvert = \displaystyle\prod\limits_{i=1}^{(n-1)} \gamma_{t}(A^{i})$.

For $p=1,\ldots,|B|$, let $Z_{p} = D \cap (B_{p} \times A^{n})$, and
\begin{align*}
S_{p} = \big\{ {D}^{n}_{x} ~\rvert~ &\text{the submatrix of $F$ determined by $B_{p}\times {D}^{n}_{x}$ is an $n$-matrix},\\
&\text{with $x \in \{1,\ldots,\gamma_{t}(A^{n})\}$}\big\}~.
\end{align*}
Note that if $q \in Q$ is a $n$-matrix, then the projection of $q$ on $A^n$
is \emph{non-self-dominated} by the projection of $D$ on $A^n$ (the same condition used in Claim \ref{clm_di_out}).  
Moreover, if $q$ is written as $B_{p}\times D^n_x$ for some $p \in \{1,\ldots,|B|\}$ and $x \in \{1,\ldots,\gamma_{t}(A^{n})\}$,
then $D^n_x$ is non-self-dominated by the projection of $Z_p$ on $A^n$.

%We now claim that for $p=1,\ldots,|B|$, $|S_{p}| \leq |Z_{p}|$. We prove this claim in a manner very similar to the proof of Claim \ref{clm_di_out}. Let
%$S _{p} = \{{D}^{n}_{i_{1}}, {D}^{n}_{i_{2}}, \ldots, {D}^{n}_{i_{t}}\}$
%and let $\varPhi_{A^n}(Z_p)$ be the projection of $Z_{p}$ on $A^n$. As in Claim \ref{clm_di_out}, $\varPhi_{A^n}(Z_p)$ dominates
%$\cup _{x=1}^{t}{D}^{n}_{i _{x}}$, and as before,
%$\varPhi_{A^n}(Z_p) \cup \big\{ {u}^{n} _{i} ~\rvert~ i\notin \{i _{1}, i _{2}, \ldots, i _{t} \} \big \}$ is a total dominating set of $A^{n}$.
%Therefore, $ |\varPhi_{A^n}(Z_p) \cup \big\{ {u}^{n} _{i} ~\rvert~ i\notin \{i _{1}, i _{2}, \ldots, i _{t}\} \big\}| =
%|\varPhi_{A^n}(Z_p)| + ( \gamma_{t}(A^{n}) - t ) \geq \gamma_{t}(A^{n}) $.
%Hence, $t = |S_p| \leq |\varPhi_{A^n}(Z_p)| \leq | Z _{p} |$.

We now claim that for $p=1,\ldots,|B|$, $|S_{p}| \leq |Z_{p}|$.
We prove this claim in a manner very similar to the proof of Claim \ref{clm_di_out}.
Let $S _{p} = \{{D}^{n}_{i_{1}}, {D}^{n}_{i_{2}}, \ldots, {D}^{n}_{i_{t}}\}$
and let $\varPhi_{A^n}(Z_p)$ be the projection of $Z_{p}$ on $A^n$.
As in Claim \ref{clm_di_out}, $\varPhi_{A^n}(Z_p)$ dominates
$\cup _{x=1}^{t}{D}^{n}_{i _{x}}$, and if $\varPhi_{A^n}(Z_p) \cap \big\{ {u}^{n} _{i} ~\rvert~ i\notin \{i _{1}, i _{2}, \ldots, i _{t} \} \big \}$ is non-empty, let $\varPhi_{A^n}(Z_p) \cap \big\{ {u}^{n} _{i} ~\rvert~ i\notin \{i _{1}, i _{2}, \ldots, i _{t} \} \big \} = \{{u}^{n} _{i_{t+1}},\ldots, {u}^{n} _{i_{l}} \}$. Then, as before,
$\varPhi_{A^n}(Z_p) \cup \big\{ {u}^{n} _{i} ~\rvert~ i\notin \{i _{1}, i _{2},
\ldots,i_t, \ldots, i _{l} \} \big \}$
is a total dominating set of $A^{n}$, and
the sets $\varPhi_{A^n}(Z_p)$ and $\big\{ {u}^{n} _{i} ~\rvert~ i\notin \{i _{1}, i _{2}, \ldots, i _{l} \} \big \}$
are disjoint.
Therefore, $ |\varPhi_{A^n}(Z_p) \cup \big\{ {u}^{n} _{i} ~\rvert~ i\notin \{i _{1}, i _{2}, \ldots, i _{l}\} \big\}| =
|\varPhi_{A^n}(Z_p)| + ( \gamma_{t}(A^{n}) - l ) \geq \gamma_{t}(A^{n}) $.
Hence, $t = |S_p| \leq l \leq |\varPhi_{A^n}(Z_p)| \leq | Z _{p} |$.

Now, $| d_{n}| = \displaystyle\sum\limits_{p=1}^{|B|} |S_{p} | \leq \displaystyle\sum\limits_{p=1}^{|B|} |Z _{p} | \leq |D|$.
\end{proof}

To conclude the proof,
$\displaystyle\prod\limits_{i=1}^{n} \gamma_{t}(A^{i}) \leq \displaystyle\sum\limits_{j=1}^{n} | d_{j}| \leq   n| D| =  n\gamma_{t}(A^{1} \Box \cdots \Box A^{n})$.
\end{proof}

%%%%%%%%%%%%%%%%%%%%%%%%%%%%%%%%%%%%%%%%%%%%%%%%%%%%%%%%%%%%%%%%%%%%%%%%
\subsection{Proof of Theorem~\ref{thm4}}
\begin{proof}
Let ${\{x_{1},y_{1},\ldots, x_{k}, y_{k}\} }$ be a
$\gamma _{pr}$-set of $G$, where for each $i$,
$(x_{i}, y_{i}) \in E(G)$. Thus, $\gamma_{pr}(G) = 2k$.
Partition $V(G)$ into sets $D _{1} ,\ldots, D _{k}$, such that
$\{x_{i}, y_{i} \} \subseteq D_{i} \subseteq N_{G}[x_{i}, y_{i}]$
for $1 \leq i \leq k$. Similarly,
let $\{\overline{x}_{1},\overline{y}_{1},\ldots, \overline{x}_{l}, \overline{y}_{l}\}$
be a $\gamma _{pr}$-set of $H$, where for each $j$, $(\overline{x}_{j}, \overline{y}_{j}) \in E(H)$. Thus, $\gamma_{pr}(H) = 2l$.  Partition $V(H)$ into sets
$\overline{D} _{1} ,\ldots, \overline{D} _{l}$, such that
$\{ \overline{x}_{j}, \overline{y}_{j} \} \subseteq \overline{D}_{j} \subseteq
N_{H}[\overline{x}_{j}, \overline{y}_{j}]$ for $1 \leq j \leq l$.
Now, ${\{D _{1} ,\ldots, D _{k}\} } \times {\{\overline{D} _{1} ,\ldots, \overline{D} _{l}\} }$
forms a partition of $V(G\Box H)$.

Let $D$ be a $\gamma_{pr}$-set of $G \Box H$. Then, for each $gh \notin D$, either $N_{\underline{\mathbf{G}} \Box H}(gh) \cap D$
or $N_{G \Box \underline{\mathbf{H}}}(gh) \cap D$ is non-empty. Based on this observation, we define the binary $|V(G)| \times |V(H)|$ matrix $F$ such that:
\[
 F(g,h) =
 \begin{cases}
 1 & \text{if } gh\in D\text{ or } N_{G \Box \underline{\mathbf{H}}}(gh) \cap D \neq \emptyset~, \\
 0 & \text{otherwise}~.
 \end{cases}
\]

Since $D$ is a $\gamma_{pr}$-set, the subgraph of $G \Box H$ induced by $D$ has a perfect matching. Thus, $D$ can be written as the disjoint union of
\begin{align*}
D_{G} &= \{gh \in D ~\rvert~ \text{the matching edge incident to $gh$ is a \textbf{G}-edge}\}, \text{ and }\\
D_{H} &= \{gh \in D ~\rvert~ \text{the matching edge incident to $gh$ is an~\textbf{H}-edge}\}~.
\end{align*}

For $i = 1,\ldots, k$, let $Z_{G _{i}} = D_{G} \cap (D _{i} \times V(H))$,
and $Z_{H_{i}} = D_{H} \cap (D _{i} \times V(H))$. For $j = 1,\ldots, l$, let $\overline{Z}_{G_{j}} = D_{G} \cap (V(G) \times \overline{D} _{j})$, and let $\overline{Z}_{H_{j}} = D_{H} \cap (V(G) \times \overline{D} _{j})$. By Claims \ref{clm_dj} and \ref{clm_di_out}, if the submatrix of $F$ determined by $D _{i} \times \overline{D} _{j}$ 
satisfies Prop.~\ref{fact1}a, then $\overline{D} _{j}$
is dominated by $\varPhi_{H}(Z_{G_{i}}\cup Z_{H_{i}})$, and if the submatrix of $F$ determined by $D _{i} \times \overline{D} _{j}$ satisfies Prop.~\ref{fact1}b, then $D _{i}$ is
dominated by $\varPhi_{G}(\overline{Z}_{G_{j}}\cup \overline{Z}_{H_{j}})$ .

For $i= 1,\ldots, k$, and $j=1,\ldots,l$, let
\begin{align*}
S _{i} = \big\{\overline{D} _{x} ~\rvert~ &\text{the submatrix of $F$ determined by $D _{i}\times \overline{D} _{x}$ satisfies Prop.~\ref{fact1}a},\\
&\text{with $x \in \{1,\ldots, l\}$} \big\}~,\\
\overline{S} _{j} = \big\{D _{x} ~\rvert~&\text{the submatrix of $F$ determined by $D _{x}\times \overline{D} _{j}$ satisfies Prop.~\ref{fact1}b},\\
&\text{with $x \in \{ 1,\ldots, k\}$}\big\}~.
\end{align*}

Finally, let $d_H = \sum_{i=1}^{k} |S _{i} |$, and
$d_G = \sum_{j=1}^{l} |\overline{S} _{j} |$.
Then, as before, $kl \leq d_H + d_G $, since each of the $kl$ submatrices of $F$ determined by $D _{i}\times \overline{D} _{j}$ satisfies one (or both) of the conditions of Prop.~\ref{fact1}. We now prove a claim that will allow us to bound the sizes of our various sets and conclude the proof.

%\newline
%Now, we can utilize a Lemma due to Hou and Jiang~\cite{Hou:Jiang:10}
%\begin{lem}[Hou and Jiang~\cite{Hou:Jiang:10}]
%Let G and H be graphs without isolated vertices. Let $Q \subseteq V(G\Box H)$
%such that $\varPhi_{H} (Q)$ dominates H, and $Q = Q_{1}\cup Q_{2}$ where
%$Q_{1}$ has a perfect matching in $G\Box H$. Then
%$\gamma _{pr}(H) \leqslant \lvert Q_{1} \vert + 2\lvert Q_{2}\lvert$.
%\end{lem}

\begin{claim}\label{clm_pr_si_zi}For $i=1\ldots,k$,~$2| S _{i} | \leq 2|Z_{G_{i}}| + |Z_{H_{i}}|$~.
\end{claim}
\begin{proof}
Let $S _{i} = \{\overline{D} _{j _{1}}, \overline{D} _{j _{2}}, \ldots, \overline{D} _{j _{t}}\}$.
Let $A = \varPhi_{H}(Z_{G_i})$, $B = \varPhi_{H}(Z_{H_i})$, and
$C = \big\{ \overline{x}_{j} ~\rvert~ j\notin \{j _{1}, j _{2},\ldots, j _{t}\}\big\}
\cup \big\{ \overline{y} _{j} ~\rvert~ j\notin \{j _{1}, j _{2}, \ldots, j _{t}\}\big\}$.

Let $M$ be the matching on $B \cup C$ formed by taking all of the $\{\overline{x}_{j}, \overline{y}_{j}\}$ edges induced by the vertices in $C$, and then adding the edges from a maximal matching on the remaining unmatched vertices in $B$. Then, $E = A \cup B \cup C$ is a dominating set of $H$ with $M$ as a matching. Let $M_1 = V(M)$ and
$M_2 = (B\cup C)\backslash M_1$. We note that $M_1$ consists of all the vertices in $C$ plus the matched vertices from $B$, and $M_2$ contains only the unmatched vertices from $B$. Therefore, $|M_1| + 2|M_2| \leq |C| + |Z_{H_i}|$. To see this more clearly, consider a vertex $gh \in Z_{H_i}$ that is matched by an H-edge to a vertex $gh'$ such that $h \notin V(M)$. This implies that either $h'$ coincides with a vertex of $C$, or $h'$ coincides with the projection of some other vertex of $Z_{H_i}$ (because otherwise $h$ would be matched with $h'$). Therefore, $2|M_2|$ is equivalent to counting $h'$, and we see that $|M_1| + 2|M_2| \leq |C| + |Z_{H_i}|$.

In order to obtain a perfect matching of $E$, we recursively modify $E$ by choosing an unmatched vertex $h$ in $E$ (a vertex in either $A$ or $B$, since all vertices in $C$ are automatically matched), and then either matching it with an appropriate vertex, or removing it from $E$. Specifically, if $N_H(h) \backslash V(M)$ is non-empty, there exists a vertex $h' \in N_H(h) \backslash V(M)$ such that we can add $h'$ to $E$ and $(h, h')$ to the matching $M$. Otherwise, $h$ is incident on only matched vertices, and we can remove $h$ from $E$ without altering the fact that $E$ is a dominating set.

Our recursively modified $E$ (denoted by $E_{\text{rec}}$) is now a paired dominating set of $H$. Furthermore, in the worst case, we have doubled the unmatched vertices from $B$, and also doubled the vertices in $A$. Thus,
\begin{align*}
2l \leq |E_{\text{rec}}| \leq 2\lvert A\lvert + \lvert M_1\lvert + 2\lvert M_2\lvert~.
\end{align*}
Since $|M_1| + 2|M_2| \leq |C| + |Z_{H_i}|$, this implies that $2l - |C| \leq 2|A| + |Z_{H_i}|$. Furthermore, since $2l - |C| = 2|S_i|$, we see that $2|S_i| \leq 2|Z_{G_i}| + |Z_{H_i}|$~.
\end{proof}

%above is same proof written in diff way...
%\begin{proof}
%Let $A = \{u\} \times ( \{ \overline{x} _{j} \rvert j\notin \{j _{1}, j _{2}, ..., j _{t}\} \}$
%$\cup \{ \overline{y} _{j} \rvert j\notin \{j _{1}, j _{2}, ..., j _{t}\} \} )$
%where $u$ is any element of $D_{i}$.
%Then, $\varPhi_{H} (A\cup Z_{G_{i}}\cup Z_{H_{i}})$ is a dominating set of H.
%\newline
%Let $B =  A\backslash (Z_{G_{i}}\cup Z_{H_{i}})$.
%Note that subgraph induced by A has a perfect matching.
%Let C be set of those vertices of B which have their matching vertices in B,
%and $D = (B\backslash C)$.
%($\lvert A\lvert \geqslant \lvert C\lvert + 2\lvert D\lvert$).
%Let $Q_{1} = Z_{H_{i}}\cup C$ and, $Q_{2} = Z_{G_{i}} \cup D$.
%As $A\cup Z_{G_{i}}\cup Z_{H_{i}} = Q_{1}\cup Q_{2}$, therefore,
%$\varPhi_{H} (Q_{1}\cup Q_{2})$ is a dominating set of H where $Q_{1}$
%has a perfect matching in $D_{i} \times H$.
%\newline
%Now by above lemma, $\lvert Q_{1}\lvert + 2\lvert Q_{2}\lvert$
%$\geqslant \gamma _{pr}(H)$.
%\newline
%But, $\lvert Q_{1}\lvert + 2\lvert Q_{2}\lvert$
%= $\lvert Z_{H_{i}}\lvert + 2\lvert Z_{G_{i}}\lvert + \lvert C\lvert$ +
%$2\lvert D\lvert$.
%\newline
%Therefore, $\lvert Z_{H_{i}}\lvert + 2\lvert Z_{G_{i}}\lvert + \lvert A\lvert$
%$\geqslant \gamma _{pr}(H)$.
%\newline
%As $\lvert A\lvert = \gamma _{pr}(H) - 2t$,
%$\lvert Z_{H_{i}}\lvert + 2\lvert Z_{G_{i}}\lvert \geqslant \gamma _{pr}(H) - \lvert A\lvert$
%$\geqslant \gamma _{pr}(H) - (\gamma _{pr}(H) - 2t)$.
%\newline
%Hence, $\lvert Z_{H_{i}}\lvert + 2\lvert Z_{G_{i}}\lvert \geqslant 2\lvert S _{i} \lvert$
%\end{proof}
Similarly, for $j=1,\ldots,l$,~we can show that $2|\overline{S}_{j}| \leq |\overline{Z}_{G_{j}}| + 2|\overline{Z}_{H_{j}}|$~. We now see
\begin{align*}
 2 \displaystyle\sum\limits_{i=1}^{k} | S _{i} | +
 2\displaystyle\sum\limits_{j=1}^{l} | \overline{S} _{j} |  &\leq 2\displaystyle\sum\limits_{i=1}^{k} |Z_{G_{i}} |
 + \displaystyle\sum\limits_{i=1}^{k} |Z_{H_{i}} |
 + \displaystyle\sum\limits_{j=1}^{l} |\overline{Z}_{G_{j}} |
 + 2\displaystyle\sum\limits_{j=1}^{l} | \overline{Z}_{H_{j}} |~,\\
&\leq \underbrace{\displaystyle\sum\limits_{i=1}^{k} |Z_{G_{i}} |
 + \displaystyle\sum\limits_{i=1}^{k} |Z_{H_{i}} |}_{D}
 + \underbrace{\displaystyle\sum\limits_{j=1}^{l} |\overline{Z}_{G_{j}} |
 + \displaystyle\sum\limits_{j=1}^{l} | \overline{Z}_{H_{j}} |}_{D} + 
\underbrace{\displaystyle\sum\limits_{i=1}^{k} |Z_{G_{i}} | + \displaystyle\sum\limits_{j=1}^{l} | \overline{Z}_{H_{j}} |}_{D}~,\\
&\leq 3|D|~.
\end{align*}
To conclude the proof, we note that 
\begin{align*}
2( d_H + d_G ) = 2 \displaystyle\sum\limits_{i=1}^{k} | S _{i} | +
 2\displaystyle\sum\limits_{j=1}^{l} | \overline{S} _{j} | &\leq 3|D|~,\\
2(kl) = \gamma _{pr}(G) \frac{\gamma_{pr}(H)}{2} &\leq 3|D|~,\\
\gamma _{pr}(G) \gamma_{pr}(H) \leq 6 \gamma_{pr}(G\Box H)~.
\end{align*}
\end{proof}

%%%%%%%%%%%%%%%%%%%%%%%%%%%%%%%%%%%%%%%%%%%%%%%%%%%%%%%%%%%%%%%%%%%%%%%%
\subsection{Proof of Theorem~\ref{thm_pr_n}}
\begin{proof}
For $i=1,\ldots,n$, let $k_i = \gamma_{pr}(A^{i})/2$, and let $\{x^{i}_{1},y^i_1 ,\ldots, x^{i}_{k_i}, y^{i}_{k_i} \}$
be a $\gamma _{pr}$-set of $A^{i}$, and $D^{i}_1 ,\ldots, D^{i} _{k_i}$
be the corresponding partitions (as defined in Theorem \ref{thm4}).

Let $Q = \{D^{1} _{1} ,\ldots, D^{1} _{k_1}\} \times \cdots \times \{D^{n} _{1} ,\ldots, D^{n} _{k_n}\}$.
Then $Q$ forms a partition of $V(A^{1} \Box \cdots \Box A^{n})$ with 
$|Q| = \displaystyle\prod\limits_{i=1}^{n} \gamma_{pr}(A^{i})/2 =  \frac{1}{2^n}\displaystyle\prod\limits_{i=1}^{n} \gamma_{pr}(A^{i})$. 

Let $D$ be a $\gamma_{pr}$-set of $A^{1} \Box \cdots \Box A^{n}$. Then, for each $u\in V(A^{1} \Box \cdots \Box A^{n})$, there exists an $i$ such that $N _{\Box A^i}(u) \cap D$ is non-empty. We now proceed slightly differently than previously. Based on this observation (as in the 2-dimensional case), we define $n$ different matrices $F^i$ with $i = 1,\ldots, n$, where each of the $n$ matrices is an $n$-ary $|V(A^1)| \times \cdots \times |V(A^n)|$ matrix $F^i$ such that:
\[
 F^i(u_1,\ldots,u_n) =
 \begin{cases}
 i & \text{if } u_1\cdots u_n \in D~,\\
j_{\min} & \text{where~} j_{\min} = \min \{~j~\lvert~N _{\Box A^j}(u_1\cdots u_n) \cap D \neq \emptyset \}~.
 \end{cases}
\]
Thus, each of the $n$ matrices $F^i$ with $i = 1,\ldots,n$ differs only in the entries that correspond to vertices in the paired dominating set $D$.

For $j=1,\ldots,n$ and $i = 1,\ldots,n$,  let $d^i_{j} \subseteq Q$ be
the set of the elements in $Q$ which are $j$-matrices in the matrix $F^i$.
By Prop.~\ref{fact2}, each element of $Q$ belongs to at least one $d^i_{j}$-set
for each $i = 1,\ldots,n$.
Now, if an element $q \in Q$ belongs to the $d^i_{j}$-set, then
$q$ also belongs to the $d^j_{j}$-set. To see this, if $M_i$ and $M_j$ are the submatrices determined by $q$ with respect to the matrices $F^i$ and $F^j$, respectively, then all the entries that do not
match in $M_i$ and $M_j$ have value $j$ in $M_j$.
Thus, each $q \in Q$ belongs to at least one $d^i_{i}$-set
for some $i \in \{1,\ldots,n\}$.
Then, $\displaystyle \frac{1}{2^n} \prod\limits_{i=1}^{n} \gamma_{pr}(A^{i}) \leq
\displaystyle\sum\limits_{i=1}^{n} | d^i_{i}|$.

Similar to $Q$, let $B = \{D^{1} _{1}, \ldots, D^{1} _{k_1}\} \times \cdots \times \{D^{n-1} _{1} ,\ldots, D^{n-1} _{k_{n-1}}\}$.
For convenience, we denote $B$ as $\{B_{1}, \ldots, B_{|B|}\}$, where $\rvert B\rvert = \displaystyle\prod\limits_{i=1}^{n-1} \gamma_{pr}(A^{i})/2 =\frac{1}{2^{n-1}}\displaystyle\prod\limits_{i=1}^{n-1} \gamma_{pr}(A^{i})$.

Since $D$ is a $\gamma_{pr}$-set, the subgraph of $A^{1} \Box \cdots \Box A^{n}$ induced by $D$ has a perfect matching. Let
\begin{align*}
D_{i} &= \{u \in D ~\rvert~ \text{the matching edge incident to $u$ is in $E_i$}\}~.
\end{align*}
Then, $D$ can be written as the disjoint union of the subsets $D_i$. For $p=1,\ldots,|B|$ and $i = 1,\ldots,n$, let $Z^i_{p} = D_i \cap (B_{p} \times A^{n})$, and
\begin{align*}
%S^i_{p} &= \big\{ {D}^{n}_{x} ~\rvert~ \text{the submatrix of $F^n$ determined by $B_{p}\times {D}^{n}_{x}$ is an $i$-matrix, with $x \in \{1,\ldots,k_n\}$}\big\}~.
S_{p} = \big\{ {D}^{n}_{x} ~\rvert~ &\text{the submatrix of $F^n$ determined by $B_{p}\times {D}^{n}_{x}$ is an $n$-matrix},\\
&\text{with $x \in \{1,\ldots,k_n\}$}\big\}~.
\end{align*}

\begin{claim}\label{claim_n_pr_dj_d}
For $p=1,\ldots,|B|$, 
$2|S_p| \leq 2\lvert Z^1_p\lvert + \cdots + 2\lvert Z^{n-1}_p\lvert + |Z^n_{p}|$. 
\end{claim}
\begin{proof}Let
$S_{p} = \{{D}^{n}_{j_{1}}, {D}^{n}_{j_{2}}, \ldots, {D}^{n}_{j_{t}}\}$, and for $j = 1,\ldots, n$, let $V_j = \varPhi_{A^n}(Z^j_{p})$. Note that $|V_j| \leq |Z^j_{p}|$. Similiar to the proof of Claim \ref{clm_pr_si_zi}, let $C = \big\{ x^n_{j} ~\rvert~ j\notin \{j _{1}, j _{2},\ldots, j _{t}\}\big\}
\cup \big\{ y^n_{j} ~\rvert~ j\notin \{j _{1}, j _{2}, \ldots, j _{t}\}\big\}$.

Let $M$ be the matching on $V_n \cup C$ formed by taking all of the $\{x^n_{j}, y^n_{j}\}$ edges induced by the vertices in $C$, and then adding the edges from a maximal matching on the remaining unmatched vertices in $V_n$. Then, $E = V_1 \cup \cdots \cup V_n \cup C$ is a dominating set of $A^n$ with $M$ as a matching.

Let $M_1 = V(M)$ and
$M_2 = (V_n \cup C)\backslash M_1$. We note that $M_1$ consists of all the vertices in $C$ plus the matched vertices from $V_n$, and $M_2$ contains only the unmatched vertices from $V_n$.

In order to obtain a perfect matching, we recursively modify $E$ by choosing an unmatched vertex $a$ in $E$, and then either matching it with an appropriate vertex, or removing it from $E$. Specifically, if $N_{A^n}(a) \backslash V(M)$ is non-empty, there exists a vertex $a' \in N_{A^n}(a) \backslash V(M)$ such that we can add $a'$ to $E$ and $(a, a')$ to the matching $M$. Otherwise, $a$ is incident on only matched vertices, and we can safely remove it from $E$ without altering the fact that $E$ is a dominating set.

Our recursively modified $E$ (denoted by $E_{\text{rec}}$) is now a paired dominating set of $A_n$. Furthermore, in the worst case, we have doubled the unmatched vertices from $V_n$, and also doubled the vertices in $V_1,\ldots,V_{n-1}$. Thus,
\begin{align*}
2k_n \leq |E_{\text{rec}}| \leq 2\lvert V_1\lvert + \cdots + 2\lvert V_{n-1}\lvert+ \lvert M_1\lvert + 2\lvert M_2\lvert~.
\end{align*}
This implies that $2k_n - |C| \leq 2\lvert V_1\lvert + \cdots + 2\lvert V_{n-1}\lvert + |Z^n_{p}|$.
Since $2 k_n - |C| = 2|S_p|$, therefore,
$2|S_p| \leq 2\lvert V_1\lvert + \cdots + 2\lvert V_{n-1}\lvert + |Z^n_{p}|
\leq 2\lvert Z^1_p\lvert + \cdots + 2\lvert Z^{n-1}_p\lvert + |Z^n_{p}|$~.
\end{proof}

To conclude the proof, we follow a similar method as in the proof of Theorem 4. We begin by noting that,
\begin{align*}
|d^n_{n}| &= \displaystyle\sum\limits_{p=1}^{|B|} |S_{p} |~.
\end{align*}

Using Claim \ref{claim_n_pr_dj_d}, we now see
\begin{align*}
2\displaystyle\sum\limits_{p=1}^{|B|} |S_{p} | & \leq
\displaystyle\sum\limits_{p=1}^{|B|} \Big(2\sum\limits_{j=1}^{n} |Z^j_{p}| - |Z^n_{p}|\Big)
= 2|D| - \displaystyle\sum\limits_{p=1}^{|B|} |Z^n_{p}| = 2|D| - |D_{n}|~.
\end{align*}

Therefore, $2|d^n_{n}| \leq 2|D| - |D_{n}|$.
Similarly, we can show that $2|d^i_{i}| \leq 2|D| - |D_{i}|$ for $i = 1, \ldots, n$. To conclude the proof, we see
\begin{align*}
\frac{1}{2^{n-1}}\displaystyle\prod\limits_{i=1}^{n} \gamma _{pr}(A_i) = 2(k_1\cdots k_n)
&\leq 2\displaystyle\sum\limits_{i=1}^{n}| d^i_{i}| \leq
2n|D| - \displaystyle\sum\limits_{i=1}^{n} |D_{i}| = (2n - 1)|D|~,\\
\displaystyle\prod\limits_{i=1}^{n} \gamma _{pr}(A_i) &\leq 2^{n-1}(2n - 1)\gamma _{pr}(A_1 \Box \cdots \Box A_n)~.
\end{align*}

\end{proof}

%%%%%%%%%%%%%%%%%%%%%%%%%%%%%%%%%%%%%%%%%%%%%%%%%%
\section*{Acknowledgements} The authors would like to acknowledge the support of NSF-CMMI-0926618, the Rice University VIGRE program (NSF DMS-0739420 and EMSW21-VIGRE), and the Global Initiatives Fund (Brown School of Engineering at Rice University), under the aegis of SURGE (Summer Undergraduate Research Grant for Excellence), a joint program with the ITT Kanpur and the Rice Center for Engineering Leadership.

%Keerti Choudhary was supported by the Global Initiatives Fund at 
%Rice's Brown School of Engineering, under the aegis of SURGE,
%an undergraduate research program at IIT Kanpur and RCEL,
%the Rice Center for Engineering Leadership.

% code for adding to source control
%svn import KT https://cvs.caam.rice.edu:3129/svn/hicks_workgroup/Papers/Keerti/trunk -m "First Import"
%svn checkout https://cvs.caam.rice.edu:3129/svn/hicks_workgroup/Papers/Keerti/trunk Keerti
%svn propset svn:eol-style native *
%>  sm19
%>  h7GZtv3^

%\bibliography{b1}
\bibliographystyle{plain}
\bibliography{vizing}

\begin{thebibliography}{1}

\bibitem{viz_survey_2009}
B.Stjan, B.Bre\u{s}ar, P.Dorbec, W.Goddard, B.Hartnell, M.Henning,
  S.Klav\u{z}ar, and D.Rall.
\newblock Vizing's conjecture: A survey and recent results, 2009.
\newblock preprint.

\bibitem{viz_clark_suen}
W.~Clark and S.~Suen.
\newblock An inequality related to {V}izing's conjecture.
\newblock {\em Electronic Journal of Combinatorics}, 7(Note 4), 2000.

\bibitem{ho_total_dom}
P.~T. Ho.
\newblock A note on the total domination number.
\newblock {\em Utilitas Mathematica}, 77:97--100, 2008.

\bibitem{hou_jiang_pr_dom}
X.~M. Hou and F.~Jiang.
\newblock Paired domination of cartesian products of graphs.
\newblock {\em Journal of Mathematical Research \& Exposition}, 30(1):181--185,
  2010.

\end{thebibliography}

\end{document}